\documentclass[12pt,leqno]{amsart}
\usepackage{amsmath}
\usepackage{amssymb}
\def\overset#1#2{{\mathrel{\mathop {{#2}_{}}\limits^{#1}}}}
\def\underset#1#2{{\mathrel{\mathop {{}_{} {#2}}\limits_{{#1}_{}}}}}
\def\upplim_#1{\underset{#1}{\overline\lim}\;}
\def\lowlim_#1{\underset{#1}{\underline\lim}\;}
\newtheorem{corollary}[equation]{Corollary}

\newtheorem{lemma}[equation]{Lemma}
\newtheorem{proposition}[equation]{Proposition}

\newtheorem{theorem}[equation]{Theorem}

\newcommand{\B}{\mathbb{B}}
\newcommand{\C}{{\mathbb{C}}}

\renewcommand{\P}{{\mathbb{P}}}

\newcommand{\ddc}{\mathrm{dd^c}}

\newcommand{\ric}{\mathrm{Ric}}
\newcommand{\supp}{\mathrm{Supp}\,}

\numberwithin{equation}{section}

\title[Non-integrated defect relation for meromorphic maps]{Improvement of non-integrated defect relation for meromorphic maps from K\"{a}hler manifolds} 

\begin{document}
\author{Tran Duc Ngoc}
\author{Si Duc Quang}

\def\thefootnote{\empty}
\footnotetext{
\hskip8pt This research is funded by Vietnam National Foundation for Science and Technology Development (NAFOSTED) under grant number 101.02-2021.12.}

\begin{abstract}
The purpose of this paper is to establish a non-integrated defect relation for meromorphic mappings from a complete K\"{a}hler manifold into a projective variety intersecting an arbitrary family of hypersurfaces with explicit truncation level. In our result, both the total defect and the truncation level are estimated independently of the number of involving hypersurfaces. Our result generalizes and improves the previous results in this topic. 
\end{abstract}

\maketitle

\section{Introduction and Main result} 

Let $V$ be a $k$-dimension subvariety of $\P^n(\C)$. Let $\mathcal D=\{D_1,\ldots,D_q\}$ be a family of hypersurfaces in $\P^n(\C)$. We say that the family $\mathcal D$ is in $N$-subgeneral position with respect to $V$ if
$$ V\cap\bigcap_{j=0}^ND_{i_j}\ \forall 1\le i_0<\cdots<i_{N}\le q.$$
If $\mathcal D$ is in $k$-subgeneral position with respect to $V$ then we say that it is in general position with respect to $V$.
 
Let $M$ be a complete K\"{a}hler manifold of dimension $m$. Let $f: M\longrightarrow V$ be a meromorphic mapping. We denote by and $\Omega_f$ the pull-back of the Fubini-Study form $\Omega$ on $\P^n(\C)$ by $f$. For a positive integer $\mu_0$ and a hypersurface $D$ of degree $d$ in $\P^n(\C)$ with $f(M)\not\subset D$. Denote by $\nu_f(D)(p)$ the intersection multiplicity of the image of $f$ and $D$ at $f(p)$ for each point $p\in M$.

In 1985, H. Fujimoto \cite{F85} defined the notion of non-integrated defect of $f$ with respect to $D$ truncated to level $\mu_0$ by
$$\delta_f^{[\mu_0]}(D):= 1- \inf\{\eta\ge 0: \eta \text{ satisfies condition }(*)\}.$$
Here, the condition (*) means that there exists a bounded non-negative continuous function $h$ on $M$ whose order of each zero is not less than $\min \{\nu_f(D), \mu_0\}$ such that
$$d\eta\Omega_f +\dfrac{\sqrt{-1}}{2\pi}\partial\bar\partial\log h^2\ge [\min\{\nu_f(D), \mu_0\}].$$
With this definition, Fujimoto gave a result analogous to the defect relation in Nevanlinna theory as follows. 

\vskip0.2cm
\noindent
{\bf Theorem  A} (see \cite[Theorem 1.1]{F85}). \ {\it Let $M$ be an $m$-dimensional complete K\"ahler manifold  and $\omega$ be a K\"{a}hler form of $M.$  Assume that the universal covering of $M$ is biholomorphic to a ball in $\mathbb C^m.$ Let $f : M\to \P^n(\C)$ be a linearly nondegenerate meromorphic map. Let $H_1,\ldots,H_q$ be hyperplanes of $\P^n(\C)$ in general position. Assume that there eixst $\rho\ge 0$ and a bounded continuous function $h\ge 0$ on $M$ such that $\rho \Omega_f +\ddc \log h^2\ge \ric \ \omega.$ Then we have
$$\sum_{i=1}^{q} \delta_{f}^{[n]} (H_i) \le n+1+ \rho n(n+1).$$}

The inequality in Theorem A is called the non-integrated defect relation with truncated multiplicity (to level $n$) of the meromorphic mapping $f$ on K\"{a}hler manifold $M$ with the family of hyperplanes $H_1,\ldots,H_q$. The quantity in the right hand side of the inequality is called the total defect. Recently, motivated by the establishment of the second main theorem of value distribution theory for hypersurface targets, the study of non-integrated defect relation of meromorphic mappings from complete K\"{a}hler manifolds into projective varieties has been much attended by many authors. We refer readers to the articles \cite{CH,Q16,Q21,RS,TT,TQ,Y13} and references therein for the development of related subjects. We state here some recent important results.

Firstly, consider the case where $f:M\rightarrow\P^n(\C)$ is an algebraically nondegenerate meromorphic mapping and $\mathcal D=\{D_1,\ldots,D_q\}$ is a family of hypersurfaces with $D_i\not\supset f(M)\ (1\le i\le q)$. By using the filtration method of P. Corvaja and U. Zannier \cite{CZ}, some authors obtained the following results:

$\bullet$ M. Ru and S. Sogome \cite{RS}; if $\mathcal D$ is in general position in $\P^n(\C)$ and $\epsilon>0$, then
$$\sum_{i=1}^q \delta^{[L-1]}_{f}(D_i)\le n+1+\epsilon+\frac{\rho L(L-1)}{d}$$
where $L=\left[2^{n^2+4n}e^nd^{2n}\left(n I\left(\epsilon^{-1}\right)\right)^n\right]$, $d=lcm(\deg D_1,\ldots,\deg D_q)$ and $[x]$ (resp. $I(x)$) denotes the smallest integer not exceed (resp. exceed) the real number $x$.

$\bullet$ Q. Yan \cite{Y13}; if $\mathcal D$ is in $N$-subgeneral position in $\P^n(\C)$ and $\epsilon>0$, then
$$\sum_{i=1}^q \delta^{[L-1]}_{f}(D_i)\le N(n+1)+\epsilon+\frac{\rho L(L-1)}{d},$$
 where $L=\left[(3eNdI(\epsilon^{-1}))^n(n+1)^{3n}\right].$

$\bullet$ S. D. Quang and N. T. Q. Phuong and N. T. Nhung \cite{Q16}; if $\mathcal D$ is in $N$-subgeneral position in $\P^n(\C)$ and $\epsilon>0$,  then
$$\sum_{j=1}^q \delta_f^{[L-1]}(D_j)\leq (N-n+1)(n+1)+\epsilon+\frac{\rho L(L-1)}{d},$$
 where $L=\left[e^{n+2}\left(d(N+n-1)(n+1)^2 I(\epsilon^{-1})\right)^n\right].$

In the above three results, the right hand side of the inequalities do not depend on the number $q$ of hypersurfaces. However, they can not be applied to the case of degenerate mappings, i.e., the case where the image of the mappings are contained in projective subvarieties of $\P^n(\C)$.

In another direction, by using the Nochka's weighs for hypersurfaces constructed by the second author and D. P. An in \cite{QA}, D. D. Thai-S. D. Quang \cite{TQ} and W. Chen-Q. Han \cite{CH} initially considered the case where the meromorphic mapping $f$ from $M$ into a $k$-dimension subvariety $V\subset\P^n(\C)$ and independently obtained: 
$$\sum_{j=1}^q \delta^{[H_V(d)-1]}_f\left(D_j\right) \leq \frac{(2 N-k+1) H_V(d)}{k+1}+\frac{\rho(2 N-k+1) H_V(d)(H_V(d)-1)}{(k+1)d},$$
where $\mathcal D$ is in $N$-subgeneral position with respect to $V$, the map $f$ is asummed to be nondegenerate over $I_d(V)$ (i.e., there is no non-trivial homogeneous polynomial $Q$ of degree $d$ which defines a hypersurface $D$ such that $f(M)\subset D$ and $D\not\supset V$) and $H_V(d)$ is the Hilbert function of $V$. In this result, the truncation level $H_V(d)$ is much smaller than that of the above three mentioned results. However, the first quantity $ \frac{(2 N-k+1) H_V(d)}{k+1}$ of the right hand side is much larger than usual. 

Recently, by using the estimate of Chow weights due to J. Evertse and R. Ferretti \cite{EF01,EF02}, S. D. Quang \textit{et all.} \cite{Q21,Q22a} considered the case of algebraically nondegenerate meromorphic mapping $f$ from $M$ into a $k$-dimension subvariety $V\subset\P^n(\C)$ and obtained the following:

$\bullet$ S. D. Quang and L. N. Quynh and N. T. Nhung \cite{Q21};  if $\mathcal D$ is in $N$-subgeneral position with respect to $V$, $\epsilon >0$, then
$$\sum_{j=1}^q \delta_f^{[L-1]}(D_j)\leq (N-k+1)(k+1)+\epsilon+\frac{\rho\epsilon L(L-1)}{d},$$ 
where $L=\left[d^{k^2+k}\deg(V)^{k+1}e^kp^k(2k+4)^k\ell^k\epsilon^{-k}+1\right]$ and $\ell=(k+1)(q!)$. 

$\bullet$ S. D. Quang \cite{Q22a};  if $\mathcal D$ has distributive constant $\Delta_{\mathcal D,V}$ with respect to $V$, $\epsilon >0$, then
$$\sum_{j=1}^q \delta_f^{[L-1]}(D_j) \leq\Delta_{\mathcal D,V}(k+1)+\epsilon+\frac{\rho\epsilon (L-1)}{d^{k+1}(2k+1)(k+1)(q!)\deg(V)},$$ 
where  $L=\left[d^{k^2+k}\deg(V)^{k+1}e^k \Delta_{\mathcal D,V}^k(2k+4)^k\ell^k\epsilon^{-k}+1\right]$ and $\ell=(k+1)(q!)$, the distributive constant $\Delta_{\mathcal D,V}$ of the family $\mathcal D$ with respect to $V$ is defined by (cf. \cite[Definition 3.3]{Q22a})
$$ \Delta_{\mathcal D,V}:=\underset{\emptyset\ne\Gamma\subset\{1,\ldots,q\}}\max\dfrac{\sharp\Gamma}{\dim V-\dim\left (V\cap\bigcap_{j\in\Gamma} D_j\right )}$$
(note that $\dim\varnothing =-\infty$). For $V=\P^n(\C)$, we will write $\Delta_{\mathcal D}$ for $\Delta_{\mathcal D,\P^n(\C)}$.

The above two results can be applied to the case of degenerate meromorphic mappings and also the quantity $\Delta_{\mathcal D,V}(k+1)$ is quite small. However,  in the coefficient of $\rho$, there is the term $(q!\epsilon^{-1})^{k-1}> ((q-1)!)^{k-1}$. Therefore, if $\rho$ is large enough then the right hand side of the defect relation always exceeds $q$, and hence the defect relations are meaningless. This term comes from the appearance of the error term in the estimate of the Hilbert weights (the last term in the inequality of Theorem \ref{2.1} when this theorem is applied).

Our aim in this paper is to improve the last above mentioned result (which is the most general case). We will reduce the above bound of the total defect (i.e, the quantity in the right hand side) so that it do not depend on $q$. In order to do so, we will apply the new estimate of Chow weights in \cite{Q22c} (see Lemma \ref{2.4}) and also give some new techniques to control the error terms occuring when the theorem on the estimate of Hilbert weights is applied. Our main result is stated as follows.
 
\begin{theorem}\label{1.1} 
Let $M$ be an $m$-dimensional complete K\"ahler manifold  and $\omega$ be a K\"ahler form of $M.$  Assume that the universal covering of $M$ is biholomorphic to a ball in $\mathbb C^m.$ Let $f$ be an algebraically nondegenerate meromorphic map of $M$ into a subvariety $V$ of dimension $k$ in $\P^n(\C)$. Let $\mathcal D=\{D_1,\ldots,D_q\}$ be a family of $q$ hypersurfaces in $\P^n(\C)$ with the distributive constant $\Delta_{\mathcal D,V}$ with respect to $V$.  Let $d = lcm\{\deg D_1,\ldots,\deg D_q\}$. Assume that for some $\rho \ge 0,$ there exists a bounded continuous function $h \geq 0$ on $M$ such that
$$\rho\Omega_f + \ddc \log h^2 \geq \ric \ \omega.$$
Then, for each  $\epsilon>0,$ we have
$$\sum_{j=1}^q {\delta}^{[L-1]}_{f} (D_j) \le \Delta_{\mathcal D,V}(k+1)+ \epsilon +\dfrac{\rho(L-1)(k+1)}{ud},$$
where $u=\lceil\Delta_{\mathcal D,V}(2k+1)(k+1)d^k\deg(V)(\Delta_{\mathcal D,V}(k+1)\epsilon^{-1}+1)\rceil$ and $L=d^k\deg(V)^k\binom{k+u}{k}$, which satisfies
$L\le d^{k^2+k}\deg(V)^{k+1}e^k\Delta_{\mathcal D,V}^k(2k+5)^k(\Delta_{\mathcal D,V}(k+1)\epsilon^{-1}+1)^k$.
\end{theorem}

Here, $\lceil x\rceil$ denotes the smallest integer which is not less than the real number $x$. 

\noindent
\textbf{Remark:} a) We see that the total defect obtained from Theorem \ref{1.1} is bounded above by the quatity:
$$\Delta_{\mathcal D,V}(k+1)+ \epsilon +\Delta_{\mathcal D,V}^{k-1}\frac{\rho d^{k^2-1}\deg(V)^{k}e^k(2k+5)^{k}(\Delta_{\mathcal D,V}(k+1)\epsilon^{-1}+1)^{k-1}}{2k+1}.$$
b) If $V=\P^n(\C)$ and $D_1,\ldots,D_q$ are in general position then $\deg(V)=1,k=n,\Delta_{\mathcal D,V}=1$ and the total defect obtained from Theorem \ref{1.1} is bounded above by 
$$n+1+\epsilon+\frac{\rho d^{n^2-1}e^n(2n+5)^{n}((n+1)\epsilon^{-1}+1)^{n-1}}{2n+1}\cong n+1+\epsilon +O(2^ne^nd^{n^2}(n\epsilon^{-1})^{n-1})\rho.$$ 
Note that, the bound above of the total defect obtained by M. Ru and S. Sogome \cite{RS} in this case is: $n+1+\epsilon+O(2^{2n^2+4n}e^{2n}d^{4n-1}(n\epsilon^{-1})^{2n})\rho.$

c) The map $f$ is said to be ramified over a hypersurface $D_i$ with multiplicity at least $m_i$ if either $f(M)\subset D_i$ or $\nu_f(D_i)(z)\ge m_i$ for all $z\in\supp\nu_f(D_i)$. By choosing $\epsilon=\frac{1}{2}$, from Theorem \ref{1.1} we obtain the following corollary.
\begin{corollary}\label{1.2} Let $M,\omega,V,\mathcal D=\{D_1,\ldots,D_q\}$ and $d$ be as in Theorem \ref{1.1}. Let $f$ be a meromorphic mapping from $M$ into $V$. Assume that for some $\rho \ge 0,$ there exists a bounded continuous function $h \geq 0$ on $M$ such that $\rho\Omega_f + \ddc \log h^2 \geq \ric \ \omega.$ If $f$ is ramified over each $D_i$ with multilicity at least $m_i$ such that
$$ \sum_{j=1}^q\left(1-\frac{L-1}{m_i}\right)> \Delta_{\mathcal D,V}(k+1)+\frac{1}{2} +\dfrac{\rho (L-1)(k+1)}{ud}, $$
where $u=\lceil\Delta_{\mathcal D,V}(2k+1)(k+1)d^k\deg(V)(2\Delta_{\mathcal D,V}(k+1)+1)\rceil$ and $L=d^k\deg(V)^k\binom{k+u}{k}$, then $f$ is algebraically degenerate.
\end{corollary}

\section{Basic notions and auxiliary results from Nevanlinna theory}

\noindent
{\bf A. Notation.}\ We set $\|z\| = \big(|z_1|^2 + \dots + |z_m|^2\big)^{1/2}$ for
$z = (z_1,\dots,z_m) \in \mathbb C^m$ and define
\begin{align*}
\B^m(r):= \{ z \in \mathbb C^m : \|z\| < r\}, \ S(r):= \{ z \in \mathbb C^m : \|z\| = r\}\ (0<r<\infty).
\end{align*}
Define 
$$v_{m-1}(z) := \big(\ddc \|z\|^2\big)^{m-1}\quad \quad \text{and}$$
$$\sigma_m(z):= {\rm d^c} \text{log}\|z\|^2 \land \big(\ddc \text{log}\|z\|^2\big)^{m-1}
 \text{on} \quad \mathbb C^m \setminus \{0\}.$$

For a divisor $\nu$ on  a ball $\B^m(R)$ in $\mathbb C^m$, with $0<r_0<r<R$, define
$$ N(r,r_0,\nu)=\int\limits_{r_0}^r \dfrac {n(t)}{t^{2m-1}}dt\ \text{ where }\ n(t) =
\begin{cases}
\int\limits_{|\nu|\,\cap \B^m(t)}
\nu(z) v_{m-1} & \text  { if } m \geq 2,\\
\sum\limits_{|z|\leq t} \nu (z) & \text { if }  m=1. 
\end{cases}$$
For a positive integer $M$ or $M= \infty$, we define $\nu^{[M]}(z)=\min\ \{M,\nu(z)\}$. Similarly, we define $n^{[M]}(t)$ and $N(r,r_0,\nu^{[M]})$ (denote it by $N^{[M]}(r,r_0,\nu)$).

For  a meromorphic function $\varphi$ on $\B^m(R)$, denote by $\nu_\varphi$ its divisor of zeros. Define
$$N_{\varphi}(r,r_0)=N(r,r_0,\nu_{\varphi}), \ N_{\varphi}^{[M]}(r,r_0)=N^{[M]}(r,r_0,\nu_{\varphi}).$$
For brevity, we will omit the character $^{[M]}$ if $M=\infty$.

Let $f : \mathbb \B^m(R) \longrightarrow \mathbb P^n(\mathbb C)$ be a meromorphic mapping with a reduced representation $\tilde f = (f_0 , \ldots , f_n)$. Set $\Vert \tilde f \Vert = \big(|f_0|^2 + \dots + |f_n|^2\big)^{1/2}$.
The characteristic function of $f$ is defined by 
$$ T_f(r,r_0)=\int_{r_0}^r\dfrac{dt}{t^{2m-1}}\int\limits_{\B^m(t)}\Omega_f\wedge v^{m-1}, \ (0<r_0<r<R). $$
By Jensen's formula, we will have
\begin{align*}
T_f(r,r_0)= \int\limits_{S(r)} \log\Vert\tilde f \Vert \sigma_m -
\int\limits_{S(r_0)}\log\Vert \tilde f\Vert \sigma_m +O(1), \text{ (as $r\rightarrow R$)}.
\end{align*}

Let $D$ be a hypersurface of degree $d$ in $\P^n(\C)$. The hypersurface $D$ has a defining homogeneous polynomial $Q$ of degree $d$ and of the following form 
$$ Q(x_0,\ldots,x_n)=\sum_{I\in\mathcal T_d}a_Ix^I,\ a_I\in\C, $$ 
where $\mathcal T_d=\{(i_0,\ldots,i_n)\in\mathbb Z^{n+1}_{\ge 0}; i_0+\cdots+i_n=d\}$, $x^I=x_0^{i_0}\cdots x_n^{i_n}$ for each $I=(i_0,\ldots,i_n)\in\mathcal T_d$. This means that
$$D=\{(x_0:\cdots:x_n)\in\P^n(\C);Q(x_0,\ldots,x_n)=0\}.$$
The proximity function of $f$ with respect to the hypersurface $D$, denoted by $m_f (r,r_0,D)$, is defined by
$$m_f (r,r_0,D)=\int_{S(r)}\log\dfrac{\|\tilde f\|^d}{|Q(\tilde f)|}\sigma_m-\int_{S(r_0)}\log\dfrac{\|\tilde f\|^d}{|Q(\tilde f)|}\sigma_m,$$
where $Q(\tilde f)=Q(f_0,\ldots,f_n)$. This definition is independent of the choice of the reduced representation $\tilde f$ of $f$. Note that, in this case $\nu_f(D)=\nu_{Q(\tilde f)}$, and by Jensen's formula, we have
$$N(r,r_0,\nu_f(D))=\int_{S(r)}\log |Q(\tilde f)|\sigma_m-\int_{S(r_0)}\log |Q(\tilde f)| \sigma_m.$$
Then, the first main theorem in Nevanlinna theory for is stated as follows
$$dT_f (r,r_0)=m_f (r,r_0,D) + N(r,r_0,\nu_f(D))+O(1),\ (r_0 < r < R),$$
where $O(1)$ is independent of $f$. For simplicity we will write $N^{[M]}_{D(f)}(r)$ and $N^{[M]}_{Q(f)}(r)$ for $N^{[M]}(r,r_0,\nu_f(D))$ and $N^{[M]}(r,r_0,\nu_{Q(\tilde f)}$ respectively.

If $ \lim\limits_{r\rightarrow 1}\sup\dfrac{T(r,r_0)}{\log 1/(1-r)}= \infty,$ then the Nevanlinna's defect of $f$ with respect to the hypersurface $D$ truncated to level $\ell$ is defined by
$$ \delta^{[\ell]}_{f,*}(D)=1-\lim\mathrm{sup}\dfrac{N^{[\ell]}_{D(f)}(r)}{\deg D \cdot T_f(r,r_0)}.$$
There is a fact that 
$$0\le \delta^{[\ell]}_f(D)\le\delta^{[\ell]}_{f,*}(D)\le 1. $$
(See Proposition 2.1 in \cite{RS})

\vskip0.2cm
\noindent
{\bf B. auxiliary results.} Using the argument in (\cite{F85}, Proposition 4.5), we have the following.

\begin{proposition}\label{2.1}
Let $F_0,\ldots ,F_{N}$ be meromorphic functions on the ball $\B^m(R_0)$ in $\mathbb C^m$ such that $\{F_0,\ldots ,F_{N}\}$ are  linearly independent over $\mathbb C.$ Then  there exists an admissible set  
$$\{\alpha_i=(\alpha_{i1},\ldots,\alpha_{im})\}_{i=0}^{N} \subset \mathbb Z^m_+$$
with $|\alpha_i|=\sum_{j=1}^{m}|\alpha_{ij}|\le i \ (0\le i \le N)$ such that the following are satisfied:

(i)\  $W_{\alpha_0,\ldots ,\alpha_{N}}(F_0,\ldots ,F_{N})\overset{Def}{:=}\det{({\mathcal D}^{\alpha_i}\ F_j)_{0\le i,j\le N}}\not\equiv 0.$ 

(ii) $W_{\alpha_0,\ldots ,\alpha_{N}}(hF_0,\ldots ,hF_{N})=h^{N+1}W_{\alpha_0,\ldots ,\alpha_{N}}(F_0,\ldots ,F_{N})$ for any nonzero meromorphic function $h$ on $\B^m(R_0).$
\end{proposition}

The following proposition is due to M. Ru and S. Sogome \cite{RS}.

\begin{proposition}[{see \cite{RS}, Proposition 3.3}]\label{2.2}
 Let $L_0,\ldots ,L_{N}$ be linear forms of $N+1$ variables and assume that they are linearly independent. Let $F$ be a meromorphic mapping of the ball $\B^m(R_0)\subset\C^m$ into $\P^{N}(\C)$ with a reduced representation $\tilde F=(F_0,\ldots ,F_{N})$ and let $(\alpha_0,\ldots ,\alpha_N)$ be an admissible set of $F$. Set $\ell=|\alpha_0|+\cdots +|\alpha_N|$ and take $t,p$ with $0< t\ell< p<1$. Then, for $0 < r_0 < R_0,$ there exists a positive constant $K$ such that for $r_0 < r < R < R_0$,
$$\int_{S(r)}\biggl |z^{\alpha_0+\cdots +\alpha_N}\dfrac{W_{\alpha_0,\ldots ,\alpha_N}(F_0,\ldots ,F_{N})}{L_0(\tilde F)\cdots L_{N}(\tilde F)}\biggl |^t\sigma_m\le K\biggl (\dfrac{R^{2m-1}}{R-r}T_F(R,r_0)\biggl )^q.$$
\end{proposition}
Here $z^{\alpha_i}=z_1^{\alpha_{i1}}\cdots z_m^{\alpha_{im}},$ where $\alpha_i=(\alpha_{i1},\ldots,\alpha_{im})\in\mathbb N^m_0$.

\vskip0.2cm 
\noindent

Let $X\subset\P^n(\C)$ be a projective variety of dimension $k$ and degree $\delta$. For $\textbf{a} = (a_0,\ldots,a_n)\in\mathbb Z^{n+1}_{\ge 0}$ we write ${\bf x}^{\bf a}$ for the monomial $x^{a_0}_0\cdots x^{a_n}_n$. Let $I=I_X$ be the prime ideal in $\C[x_0,\ldots,x_n]$ defining $X$. Let $\C[x_0,\ldots,x_n]_u$ denote the vector space of homogeneous polynomials in $\C[x_0,\ldots,x_n]$ of degree $u$ (including $0$). For $u= 1, 2,\ldots,$ put $I_u:=\C[x_0,\ldots,x_n]_u\cap I$ and define the Hilbert function $H_X$ of $X$ by
\begin{align*}
H_X(u):=\dim (\C[x_0,\ldots,x_n]_u/I_u).
\end{align*}
By the usual theory of Hilbert polynomials, we have
\begin{align*}
H_X(u)=\delta\cdot\dfrac{u^k}{k!}+O(u^{k-1}).
\end{align*}
Let ${\bf c}=(c_0,\ldots,c_n)$ be a tuple in $\mathbb R^{n+1}_{\ge 0}$ and let $e_X({\bf c})$ be the Chow weight of $X$ with respect to ${\bf c}$. The $u$-th Hilbert weight $S_X(u,{\bf c})$ of $X$ with respect to ${\bf c}$ is defined by
\begin{align*}
S_X(u,{\bf c}):=\max\left (\sum_{i=1}^{H_X(u)}{\bf a}_i\cdot{\bf c}\right),
\end{align*}
where the maximum is taken over all sets of monomials ${\bf x}^{{\bf a}_1},\ldots,{\bf x}^{{\bf a}_{H_X(u)}}$ whose residue classes modulo $I$ form a basis of $\C[x_0,\ldots,x_n]_u/I_u.$

The following theorem is due to J. Evertse and R. Ferretti \cite{EF01}.
\begin{theorem}[{Theorem 4.1 \cite{EF01}}]\label{2.3}
Let $X\subset\P^n(\C)$ be an algebraic variety of dimension $k$ and degree $\delta$. Let $u>\delta$ be an integer and let ${\bf c}=(c_0,\ldots,c_n)\in\mathbb R^{n+1}_{\geqslant 0}$.
Then
$$ \dfrac{1}{uH_X(u)}S_X(u,{\bf c})\ge\dfrac{1}{(k+1)\delta}e_X({\bf c})-\dfrac{(2k+1)\delta}{u}\cdot\left (\max_{i=0,\ldots,n}c_i\right). $$
\end{theorem}

The following lemma is due to the second author in \cite{Q22b}.
\begin{lemma}[{see \cite[Lemma 3.2]{Q22b}}]\label{2.4}
Let $Y$ be a projective subvariety of $\P^R(\C)$ of dimension $k\ge 1$ and degree $\delta_Y$. Let $\ell\ (\ell\ge k+1)$ be an integer and let ${\bf c}=(c_0,\ldots,c_R)$ be a tuple of non-negative reals. Let $\mathcal H=\{H_0,\ldots,H_R\}$ be a set of hyperplanes in $\P^R(\C)$ defined by $H_{i}=\{y_{i}=0\}\ (0\le i\le R)$. Let $\{i_1,\ldots, i_\ell\}$ be a subset of $\{0,\ldots,R\}$ such that:
\begin{itemize}
\item[(1)] $c_{i_\ell}=\min\{c_{i_1},\ldots,c_{i_\ell}\}$,
\item[(2)] $Y\cap\bigcap_{j=1}^{\ell-1}H_{i_j}\ne \varnothing$, 
\item[(3)] and $Y\not\subset H_{i_j}$ for all $j=1,\ldots,\ell$.
\end{itemize}
Let $\Delta_{\mathcal H,Y}$ be the distributive constant of the family $\mathcal H=\{H_{i_j}\}_{j=1}^\ell$ with respect to $Y$. Then
$$e_Y({\bf c})\ge \frac{\delta_Y}{\Delta_{\mathcal H,Y}}(c_{i_1}+\cdots+c_{i_\ell}).$$
\end{lemma}

\section{Proof of Theorem \ref{1.1}}

By using the universal covering if necessary, we may assume that $M=\mathbb{B}^m(R_0)$. Let $Q_1,\ldots,Q_q$ be defining homogeneous polynomials of $D_1,\ldots,D_q$ respectively, with $\deg Q_i=\deg D_i=d_i\ (1\le i\le q)$. Replacing $Q_i$ by $Q_i^{d/d_j} \ (j=1,\ldots, q)$ if necessary, we may assume that all hypersurfaces $Q_i\ (1\le i\le q)$ are of the same degree $d$. 
Consider the mapping $\Phi$ from $V$ into $\P^{q-1}(\C)$, which maps a point ${\bf x}=(x_0:\cdots:x_n)\in V$ into the point $\Phi({\bf x})\in\P^{q-1}(\C)$ given by
$$\Phi({\bf x})=(Q_1(x):\cdots : Q_{q}(x)),$$
where $x=(x_0,\ldots,x_n)$. We set $\tilde\Phi(x)=(Q_1(x),\ldots ,Q_{q}(x))$ and fix a reduced representation $\tilde f=(f_0,\ldots,f_n)$ of $f$.

Let $Y=\Phi(V)$. Since $V\cap\bigcap_{j=1}^{q}Q_j=\varnothing$, $\Phi$ is a finite morphism on $V$ and $Y$ is a complex projective subvariety of $\P^{q-1}(\C)$ with $\dim Y=k$ and of degree $\delta:=\deg(Y)\le d^{k}.\deg(V).$ 

For every ${\bf a} = (a_1,\ldots,a_q)\in\mathbb Z^q_{\ge 0}$ and ${\bf y} = (y_1,\ldots,y_q)$ we denote ${\bf y}^{\bf a} = y_{1}^{a_{1}}\ldots y_{q}^{a_{q}}$. Let $u$ be a positive integer. We set $\xi_u:=\binom{q+u}{u}$ and define
$$ Y_{u}:=\C[y_1,\ldots,y_q]_u/(I_{Y})_u, $$
which is a $\C$-vector space of dimension $H_{Y}(u)$. Put $n_u=H_{Y}(u)-1$ and let $v_0,\ldots,v_{n_u}$ be homogeneous polynomials in $\C[y_1,\ldots,y_q]_u$, whose equivalent classes form a basis of $Y_u$.  
Now, consider the meromorphic mapping $F$ from $\B^m(R_0)$ into $\P^{n_u}(\C)$ with the reduced representation given by
$$ \tilde F=(v_0(\tilde\Phi\circ \tilde f),\ldots,v_{n_u}(\tilde\Phi\circ \tilde f)).$$
Since $f$ is algebraically nondegenerate, $F$ is linearly nondegenerate.

Then there exists an admissible set $\alpha=(\alpha_0,\ldots,\alpha_{n_u})\in(\mathbb{Z}^m_+)^{n_u+1}$ such that
$$W^\alpha(F_0,\ldots,F_{n_u})=\det (D^{\alpha_i}(v_s(\Phi\circ \tilde f)))_{0\le i,s\le n_u}\neq 0.$$

Fix a point $z\in\B^m(R_0)$ such that $Q_i(\tilde f(z))\ne 0$ for all $i=1,\ldots,q$. We define 
$${\bf c}_z = (c_{1,z},\ldots,c_{q,z})\in\mathbb Z^{q},$$ 
where
$$c_{i,z}:=\log\dfrac{\|\tilde f(z)\|^d\|Q_i\|}{|Q_i(\tilde f)(z)|}\ge 0\text{ for } i=1,\ldots,q.$$
By the definition of the Hilbert weight, there are ${\bf a}_{0,z},\ldots,{\bf a}_{n_u,z}\in\mathbb N^{q}$ with
$$ {\bf a}_{j,z}=(a_{j,1,z},\ldots,a_{j,q,z}),$$ 
where $a_{j,i,z}\in\{1,\ldots,u\},$  such that the residue classes modulo $(I_Y)_u$ of ${\bf y}^{{\bf a}_{0,z}},\ldots,{\bf y}^{{\bf a}_{n_u,z}}$ form a basis of $\C[y_1,\ldots,y_q]_u/(I_Y)_u$ and
\begin{align*}
S_Y(u,{\bf c}_z)=\sum_{j=0}^{n_u}{\bf a}_{j,z}\cdot{\bf c}_z.
\end{align*}
We see that ${\bf y}^{{\bf a}_{j,z}}\in Y_u$ (modulo $(I_Y)_u$). Then we may write
$$ {\bf y}^{{\bf a}_{j,z}}=L_{j,z}(v_0,\ldots ,v_{n_u}), $$ 
where $L_{j,z}\ (0\le j\le n_u)$ are independent linear forms.
We have
\begin{align*}
\log\prod_{j=0}^{n_u} |L_{j,z}(\tilde F(z))|&=\log\prod_{j=0}^{n_u}\prod_{i=1}^q|Q_i(\tilde f(z))|^{a_{j,i,z}}\\
&=-S_Y(u,{\bf c}_z)+du(n_u+1)\log \|\tilde f(z)\| +O(u(n_u+1)).
\end{align*}
Therefore
\begin{align}\label{3.1}
S_Y(u,{\bf c}_z) = \log\prod_{j=0}^{n_u}\dfrac{1}{|L_{j,z}(\tilde F(z))|} + du(n_u+1)\log \|\tilde f(z)\|+O(u(n_u+1)).
\end{align}
From Theorem \ref{2.3} we have
\begin{align}\label{3.2}
\dfrac{1}{u(n_u+1)}S_Y(u,{\bf c}_z)\ge&\dfrac{1}{(k+1)\delta}e_Y({\bf c}_z)-\dfrac{(2k+1)\delta}{u}\max_{1\le i\le q}c_{i,z}
\end{align}
Combining (\ref{3.1}), (\ref{3.2}), we get
\begin{align}\label{3.3}
\begin{split}
\dfrac{1}{(k+1)\delta}e_Y({\bf c}_z)\le & \dfrac{1}{u(n_u+1)}\log\prod_{j=0}^{n_u}\dfrac{1}{|L_{j,z}(\tilde F(z))|} + d\log \|\tilde f(z)\| \\
&+\dfrac{(2k+1)\delta}{u}\sum_{i=1}^q\log\dfrac{\|\tilde f(z)\|^d\|Q_i\|}{|Q_i(\tilde f)(z)|}+O(1/u).
\end{split}
\end{align}
Suppose that $c_{1,z}\ge c_{2,z}\ge\cdots\ge c_{q,z}$ and denote by $\ell$ the smallest index such that $V\cap\bigcap_{j=1}^\ell D_j=\emptyset$. By Lemma \ref{2.2}, we have
\begin{align}\label{3.4}
\begin{split}
e_Y({\bf c}_z)&\ge\frac{\delta}{\Delta_{\mathcal D,V}}(c_{1,z}+\cdots +c_{\ell,z}) =\frac{\delta}{\Delta_{\mathcal D,V}}\left (\sum_{i=1}^q\log\dfrac{\|\tilde f(z)\|^d\|Q_i\|}{|Q_i(\tilde f)(z)|}\right)\\
&=\frac{\delta}{\Delta_{\mathcal D,V}}\left (\sum_{i=1}^q\log\dfrac{\|\tilde f(z)\|^d\|Q_i\|}{|Q_i(\tilde f)(z)|}\right)+O(1).
\end{split}
\end{align}
Then, from (\ref{3.3}) and (\ref{3.4}) we have
\begin{align}\label{3.5}
\begin{split}
\dfrac{1}{\Delta_{\mathcal D,V}}\log\prod_{i=1}^q\dfrac{\|\tilde f (z)\|^d}{|Q_i(\tilde f)(z)|}\le & \dfrac{k+1}{u(n_u+1)}\log\prod_{j=0}^{n_u}\dfrac{1}{|L_{j,z}(\tilde F(z))|}+d(k+1)\log \|\tilde f(z)\|\\
&+\dfrac{(2k+1)(k+1)\delta}{u}\sum_{i=1}^q\log\dfrac{\|\tilde f(z)\|^d\|Q_i\|}{|Q_i(\tilde f)(z)|}+O(1),
\end{split}
\end{align}
where the term $O(1)$ does not depend on $z$. 

Set $m_0=(2k+1)(k+1)\delta$ and $b=\dfrac{k+1}{u(n_u+1)}$. From (\ref{3.5}), we get
\begin{align}\label{3.6}
\begin{split}
\log \dfrac{\|\tilde f (z)\|^{\frac{1}{\Delta_{\mathcal D,V}}dq-d(k+1)-\frac{dm_0q}{u}}|W^\alpha(\tilde{F}(z))|^b}{\prod_{i=1}^q|Q_i(\tilde f)(z)|^{\frac{1}{\Delta_{\mathcal D,V}}-\frac{m_0}{u}}}\le b \log\dfrac{|W^\alpha(\tilde{F}(z))|}{\prod_{j=0}^{n_u}|L_{j,z}(\tilde F(z))|}+O(1).
\end{split}
\end{align}
Here we note that $L_{j,z}$ depends on $i$ and $z$, but the number of these linear forms is finite (at most $\xi_u$). We denote by $\mathcal L$ the set of all $L_{j,z}$ occuring in the above inequalities.

Then, from (\ref{3.6}) there exists a positive constant $K_0$ such that
\begin{align}\label{3.7}
\dfrac{\|\tilde f (z)\|^{\frac{1}{\Delta_{\mathcal D,V}}dq-d(k+1)-\frac{dm_0q}{u}}.|W^\alpha(\tilde{F}(z))|^b}{\prod_{i=1}^q|Q_i(\tilde f)(z)|^{\frac{1}{\Delta_{\mathcal D,V}}-\frac{m_0}{u}}}\le K_0^b.S_{\mathcal J}^b,
\end{align}
where $ S_{\mathcal J}=\dfrac{|W^\alpha(\tilde{F}(z))|}{\prod_{L\in\mathcal J}|L(\tilde F(z))|}$ for some $\mathcal J \subset \mathcal{L}$ with $\# \mathcal{J} = n_u+1$ so that $\{L\in\mathcal J\}$ is linearly independent.

We now estimate the divisor $\nu_{W^\alpha(\tilde F)}$. Fix a point $z\in\B^m(R_0)$ outside the indeterminacy locus of $f$. Suppose that
$$\max\{0,\nu_{Q_1(\tilde f)}(z)-n_u\}\ge\cdots\ge \max\{0,\nu_{Q_q(\tilde f)}(z)-n_u\}$$ 
and let $\ell$ be as above. We set $c_{i}=\max\{0,\nu_{Q_i(\tilde f)}(z)-n_u\}$ and 
$${\bf c}=(c_1,\ldots,c_q)\in\mathbb Z^q_{\ge 0}.$$
It is clear that $c_i=0$ for all $i\ge \ell$. Then there are 
$${\bf a}_j=(a_{j,1},\ldots,a_{j,p}),a_{j,s}\in\{1,\ldots,u\}, j=0,\ldots,n_u$$
such that ${\bf y}^{{\bf a}_0},\ldots,{\bf y}^{{\bf a}_{n_u}}$ is a basis of $Y_u$ and
$$ S_Y(m,{\bf c})=\sum_{j=0}^{n_u}{\bf a}_j\cdot{\bf c}.$$
Similarly as above, we write ${\bf y}^{{\bf a}_j}=L_j(v_0,\ldots,v_{n_u})$, where $L_0,\ldots,L_{n_u}$ are independent linear forms in variables $y_{i}\ (1\le i\le q)$. By the property of the general Wronskian, we have
$$W^{\alpha}(\tilde F)=cW^{\alpha}(L_0(\tilde F),\ldots,L_{n_u}(\tilde F)),$$
where $c$ is a nonzero constant. This yields that
$$ \nu_{W^{\alpha}(\tilde F)}(z)=\nu_{W^{\alpha}(L_1(\tilde F),\ldots,L_{H_Y(u)}(\tilde F))}(z)\ge\sum_{j=0}^{n_u}\max\{0,\nu_{L_j(\tilde F)}(z)-n_u\}$$
We also easily see that $ \nu_{L_j(\tilde F)}(z)=\sum_{i=1}^qa_{j,i}\nu_{Q_i(\tilde f)}(z),$ and hence
$$ \max\{0,\nu_{L_j(\tilde F)}(z)-n_u\}\ge\sum_{i=1}^{q}a_{j,i}c_{i}={{\bf a}_j}\cdot{\bf c}. $$
Thus, we have
\begin{align}\label{3.8}
 \nu_{W^{\alpha}(\tilde F)}(z)\ge\sum_{j=0}^{n_u}{{\bf a}_j}\cdot{\bf c}=S_Y(u,{\bf c}).
\end{align}
By Lemma \ref{2.4} we have
$$ e_Y({\bf c})\ge \frac{\delta}{\Delta_{\mathcal D,V}}\sum_{i=1}^{l}c_{i}= \frac{\delta}{\Delta_{\mathcal D,V}}\sum_{j=1}^{q}\max\{0,\nu_{Q_i(\tilde f)}(z)-n_u\}. $$
On the other hand, by Theorem \ref{2.1} we have 
\begin{align*}
 S_Y(u,{\bf c})&\ge\dfrac{u(n_u+1)}{(k+1)\delta}e_Y({\bf c})-(2k+1)\delta (n_u+1)\max_{1\le i\le q}c_{i}\\
&\ge\dfrac{u(n_u+1)}{(k+1)\Delta_{\mathcal D,V}}\sum_{i=1}^{q}\max\{0,\nu_{Q_i(\tilde f)}(z)-n_u\}\\
&-(2k+1)\delta(n_u+1)\sum_{i=1}^q\max\{0,\nu_{Q_i(\tilde f)}(z)-n_u\}.
\end{align*}
Combining this inequality and (\ref{3.8}), we have
\begin{align}\label{3.9}
b\nu_{W^{\alpha}(\tilde F)}(z)\ge&\left(\frac{1}{\Delta_{\mathcal D,V}}-\dfrac{m_0}{u}\right)\sum_{j=1}^{q}\max\{0,\nu_{Q_i(\tilde f)}(z)-n_u\}.
\end{align}
Therefore,
$$\left(\dfrac{1}{\Delta_{\mathcal D,V}}-\frac{m_0}{u}\right)\sum_{i=1}^{q}\nu_{Q_i(\tilde f)}(z)-b\nu_{W(\tilde F)}(z)\le \left(\dfrac{1}{\Delta_{\mathcal D,V}}-\frac{m_0}{u}\right)\sum_{i=1}^{q}\min\{\nu_{Q_i(\tilde f)}(z),n_u\}.$$
Seting $x=\frac{1}{\Delta_{\mathcal D,V}}-\frac{m_0}{u}$, from the above inequality we have
 \begin{align}\label{3.10}
\left(\dfrac{1}{\Delta_{\mathcal D,V}}-\frac{m_0}{u}\right)b\nu_{W(\tilde F)}(z)-x\left(\sum_{i=1}^{q}\nu_{Q_i(\tilde f)}(z)-\min\{\nu_{Q_i(\tilde f)}(z),n_u\}\right)\ge 0.
\end{align}

Assume that 
$$ \rho\Omega_f+\dfrac{\sqrt{-1}}{2\pi}\partial\bar\partial\log h^2\ge \ric\omega.$$
We now suppose that
$$\sum_{j=1}^q\delta^{[n_u]}_f(D_j)>\frac{k+1}{x}+\dfrac{\rho n_u(n_u+1)b}{d}.$$
Then, for each $j\in\{1,\ldots ,q\},$ there exist constants $\eta_j>0$ and continuous plurisubharmonic function $\tilde u_j$ such that 
$e^{\tilde u_j}|\varphi_j|\le \|\tilde f\|^{d\eta_j},$ where $\varphi_j$ is a holomorphic function with $\nu_{\varphi_j}=\min\{\nu_{Q_i(\tilde f)},n_u\}$ and
$$ q-\sum_{j=1}^q\eta_j>\frac{k+1}{x}+\dfrac{\rho n_u(n_u+1)b}{dx}.$$
Put $u_j=\tilde u_j+\log |\varphi_j|$, then $u_j$ is a plurisubharmonic and
$$ e^{u_j}\le \|\tilde f\|^{d\eta_j},\ j=1,\ldots ,q. $$
Let
$$v (z)=\log\left |(z^{\alpha_0+\cdots+\alpha_{n_u}})^b\dfrac{\left|W^{\alpha}(\tilde F(z))\right|^b}{\prod_{i=1}^{q}\left|Q_i(\tilde f)(z)\right|^x}\right |+x\sum_{j=1}^q u_j(z).$$
Therefore, we have the following current inequality
\begin{align*}
2\ddc[v]&\ge b[\nu_{W^{\alpha}(\tilde F)}]-x\sum_{j=1}^q[\nu_{Q_i(\tilde f)}]+x\sum_{j=1}^q2\ddc[u_j]\\
&\ge b[\nu_{W^{\alpha}(\tilde F)}]-x\left(\sum_{j=1}^q[\nu_{Q_i(\tilde f)}]-\sum_{j=1}^q[\min\{n_u,\nu_{Q_i(\tilde f)}\}]\right)\ge 0.
\end{align*}
This implies that $v$ is a plurisubharmonic function on $\B^m(R_0)$.

On the other hand, by the growth condition of $f$, there exists a continuous plurisubharmonic function $\omega\not\equiv\infty$ on $\B^m(R_0)$ such that
\begin{align*}
e^\omega {\rm d}V\le \|\tilde f\|^{2\rho}v_m
\end{align*}
Set
$$t=\dfrac{2\rho}{dx\left(q-\frac{k+1}{x}-\sum_{j=1}^q\eta_j\right)}>0$$ 
and 
$$\lambda (z)=\left|(z^{\alpha_0+\cdots +\alpha_{n_u}})^b\dfrac{\left|W^{\alpha}(\tilde F(z))\right|^b}{\left|Q_1(\tilde f)(z)\cdots Q_q(\tilde f)(z)\right|^x}\right|.$$ 

It is clear that $\dfrac{n_u(n_u+1)b}{2}t< 1,$ and the function $\zeta=\omega+ tv$ is plurisubharmonic on the K\"{a}hler manifold $M$. Choose a position number $\gamma$ such that $0<\dfrac{n_u(n_u+1)b}{2}t<\gamma<1.$
Then, we have
\begin{align}\label{3.11}
\begin{split}
e^\zeta dV&=e^{\omega +tv}dV\le e^{tv}\|\tilde f\|^{2\rho}v_m=|\lambda|^{t}(\prod_{j=1}^qe^{txu_j})\|\tilde f\|^{2\rho}v_m\\
&\le|\lambda|^t \|\tilde f\|^{2\rho+tx\sum_{j=1}^qd\eta_j}v_m=|\lambda|^t \|\tilde f\|^{dtx(q-\frac{k+1}{x})}v_m.
\end{split}
\end{align}

(a) We first consider the case where $R_0< \infty$ and $ \lim\limits_{r\rightarrow R_0}\sup\dfrac{T_f(r,r_0)}{\log 1/(R_0-r)}<\infty.$
It suffices for us to proof the Theorem in the case where $\mathbb{B}^m(R_0)=\mathbb{B}^m(1).$
Integrating both sides of (\ref{3.11}) over $\mathbb{B}^m(1)$ and using (\ref{3.7}), we have
\begin{align}\label{3.12}
\begin{split}
\int_{\B^m(1)}e^\zeta dV&\le \int_{\B^m(1)}|\lambda|^t \|\tilde f\|^{dtx(q-\frac{k+1}{x})}v_m.\\
&=2m\int_0^1r^{2m-1}\left (\int_{S(r)}\bigl (|\lambda| \|\tilde f\|^{dxq-d(k+1)}\bigl )^t\sigma_m\right )dr\\
&\le 2m\int_0^1r^{2m-1}\left (\int_{S(r)}\sum_{\mathcal J}\bigl |(z^{\alpha_0+\cdots +\alpha_{n_u}})K_0S_{\mathcal J}\bigl |^{bt}\sigma_m\right )dr,
\end{split}
\end{align}
where the summation is taken over all $\mathcal J\subset\mathcal L$ with $\sharp J=n_u+1$ and $\{L\in\mathcal J\}$ is linearly independent.

We note that $(\sum_{i=0}^{n_u}|\alpha_i|)bt\le \dfrac{n_u(n_u+1)b}{2}t<\gamma<1$. Then by Proposition \ref{2.2} there exists a positive constant $K_1$ such that, for every $0<r_0<r<r'<1,$ we have
\begin{align*}
\int_{S(r)}\left |(z^{\alpha_0+\cdots +\alpha_{n_u}})K_0S_{\mathcal J}(z)\right |^{bt}\sigma_m\le K_1\left (\dfrac{r'^{2m-1}}{r'-r}dT_f(r',r_0)\right )^{\gamma}.
\end{align*}
Choosing $r'=r+\dfrac{1-r}{eT_f(r,r_0)}$, we get $ T_f(r',r_0)\le 2T_f(r,r_0)$ outside a subset $E\subset [0,1]$ with $\int_E\dfrac{dr}{1-r}<+\infty$. Hence, the above inequality implies that
\begin{align*}
\int_{S(r)}\left |(z^{\alpha_1+\cdots +\alpha_{n_u}})K_0S_J(z)\right |^{bt}\sigma_m\le \dfrac{K}{(1-r)^\gamma}\left (\log\dfrac{1}{1-r}\right )^{\gamma}
\end{align*}
for all $r$ outside $E$, where $K$ is a positive constant. By choosing $K$ large enough, we may assume that the above inequality holds for all $r\in (0;1)$.
Then, the inequality (\ref{3.12}) yields that
\begin{align*}
\int_{\B^m(1)}e^\zeta dV&\le 2m\int_0^1r^{2m-1}\dfrac{K}{(1-r)^\gamma}\left (\log\dfrac{1}{1-r}\right )^{\gamma}dr< +\infty
\end{align*}
This contradicts the results of S. T. Yau  and L. Karp (see \cite{K82,Y76}). 

Hence, we must have
$$\sum_{j=1}^q\delta^{[n_u]}_{f}(D_j)\le\frac{k+1}{x}+\dfrac{\rho n_u(n_u+1)b}{d}=\frac{k+1}{x}+\dfrac{\rho n_u(k+1)}{ud}.$$
On the other hand, we have $u\ge\Delta_{\mathcal D,V}m_0(\Delta_{\mathcal D,V}(k+1)+\epsilon)\epsilon^{-1}$. This implies that $\frac{\Delta_{\mathcal D,V}m_0}{u}\le\frac{\epsilon}{\Delta_{\mathcal D,V}(k+1)+\epsilon}=1-\frac{\Delta_{\mathcal D,V}(k+1)}{\Delta_{\mathcal D,V}(k+1)+\epsilon}$. Therefore
\begin{align}\label{3.13}
\frac{k+1}{x}=\frac{k+1}{\frac{1}{\Delta_{\mathcal D,V}}-\frac{m_0}{u}}=\frac{\Delta_{\mathcal D,V}(k+1)}{1-\frac{\Delta_{\mathcal D,V}m_0}{u}}\le \Delta_{\mathcal D,V}(k+1)+\epsilon.
\end{align}
We have $n_u=H_Y(u)-1\le\delta\binom{k+u}{k}-1\le d^k\deg(V)\binom{k+u}{k}-1=L-1$. Then, we have
$$\sum_{j=1}^q\delta^{[L-1]}_{f}(D_j)\le\Delta_{\mathcal D,V}(k+1)+\epsilon+\dfrac{\rho (k+1)(L-1)}{ud}.$$
We estimate $L$ as follows. If $k=1$ then
\begin{align*}
L&=d\deg(V)(1+u)\\
&< d\deg(V) \left(\Delta_{\mathcal D,V}6d\deg(V)(2\Delta_{\mathcal D,V}\epsilon^{-1}+1)+2\right)\\ 
&< d^2\deg(V)^2 e\Delta_{\mathcal D,V}7(2\Delta_{\mathcal D,V}\epsilon^{-1})+1).
\end{align*}
Otherwise, if $k\ge 2$ then
\begin{align*}
L&\le d^k\deg (V)e^{k}\left(1+\dfrac{u}{k}\right)^{k}\\
&<d^k\deg (V)e^k\left(1+\dfrac{d^k\deg (V)\Delta_{\mathcal D,V}(2k+1)(k+1)(\Delta_{\mathcal D,V}(k+1)\epsilon^{-1}+1)+1}{k}\right)^k\\
&<d^{k^2+k}\deg(V)^{k+1}e^k\Delta_{\mathcal D,V}^k(2k+5)^k(\Delta_{\mathcal D,V}(k+1)\epsilon^{-1}+1)^k.
\end{align*} 

(b) We now consider the remaining case where $ \lim\limits_{r\rightarrow R_0}\sup\dfrac{T(r,r_0)}{\log 1/(R_0-r)}= \infty.$ 
It is clear that we only need to prove the following theorem.
\begin{theorem}\label{3.14}
With the assumption of Theorem {1.1}. Then, we have
$$(q-\Delta_{\mathcal D,V}(k+1)-\epsilon)T_f(r,r_0)\le \sum_{i=1}^{q}\dfrac{1}{d}N^{[L-1]}_{Q_i(\tilde f)}(r)+S(r),$$
where $S(r)$ is evaluated as follows:

(i) In the case $R_0<\infty,$ $$S(r)\le K(\log^+\dfrac{1}{R_0-r}+\log^+T_f(r,r_0))$$
for all $0<r_0<r<R_0$ outside a set $E\subset [0,R_0]$ with $\int_E\dfrac{dt}{R_0-t}<\infty$ and $K$ is a positive constant.

(ii) In the case  $R_0=\infty,$ $$S(r)\le K(\log r+\log^+T_f(r,r_0))$$
for all $0<r_0<r<\infty$ outside a set $E'\subset [0,\infty]$ with $\int_{E'} dt<\infty$ and $K$ is a positive constant.
\end{theorem}
	
\noindent\textit{Proof.} Repeating the above argument, we have
\begin{align*}
\int\limits_{S(r)}&\left|(z^{\alpha_0+\cdots +\alpha_{n_u}})^b\dfrac{\|\tilde f(z)\|^{xdq-d(k+1)}|W^{\alpha}(\tilde F)(z)|^b}{\prod_{i=1}^{q}|Q_i(\tilde f)(z)|^{x}}\right|^{t}\sigma_m\le K_1\left(\dfrac{R^{2m-1}}{R-r}dT_f(R,r_0)\right)^{\delta}.
\end{align*}
for every $0<r_0<r<R<R_0$. Using the concativity of the logarithmic function, we have
\begin{align*}
\begin{split}
b&\int_{S(r)}\log |(z^{\alpha_0+\cdots +\alpha_{n_u}})|\sigma_m+(xdq-d(k+1))\int_{S(r)}\log \|\tilde f\|\sigma_m\\
&+b\int_{S(r)}\log |W^{\alpha}(\tilde F)|\sigma_m-x\sum_{j=1}^q \int_{S(r)}\log |Q_j(\tilde f)|\sigma_m\le K\left(\log^+\dfrac{R}{R-r}+\log^+T_f(R,r_0)\right)
\end{split}
\end{align*}
for some positive constant $K$. By Jensen's formula, this inequality implies that
\begin{align}\label{3.15}
\begin{split}
\bigl(xdq&-d(k+1)\bigl)T_f(r,r_0)+bN_{W^{\alpha}(\tilde F)}(r)-x\sum_{i=1}^qN_{Q_i(\tilde f)}(r)\\
&\le K\bigl(\log^+\dfrac{R}{R-r}+\log^+T_f(R,r_0)\bigl)+O(1).
\end{split}
\end{align}
From (\ref{3.9}), we have 
$$x\sum_{i=1}^qN_{Q_i(\tilde f)}(r)-bN_{W^{\alpha}(\tilde F)}(r)\le x\sum_{i=1}^qN^{[L-1]}_{Q_i(\tilde f)}(r).$$
Combining this estimate and (\ref{3.15}), we get
\begin{align}\label{3.16}
\bigl(q-\frac{k+1}{x}\bigl)T_f(r,r_0)\le\sum_{i=1}^{q}\dfrac{1}{d}N^{[L-1]}_{Q_i(\tilde f)}(r)+\frac{K}{x}\bigl(\log^+\dfrac{R}{R-r}+\log^+T_f(R,r_0)\bigl)+O(1).
\end{align}
Choosing $R=r+\dfrac{R_0-r}{eT_f(r,r_0)}$ if $R_0<\infty$ and $R=r+\dfrac{1}{T_f(r,r_0)}$ if $R_0=\infty$, we see that
$$ T_f\bigl(r+\dfrac{R_0-r}{eT_f(r,r_0)},r_0\bigl)\le 2T_f(r,r_0)$$
outside a subset $E\subset [0,R_0)$ with $\int_E\dfrac{dr}{R_0-r}<+\infty$ in the case $R_0<\infty$ and  
$$ T_f\bigl(r+\dfrac{1}{T_f(r,r_0)},r_0\bigl)\le 2T_f(r,r_0)$$
outside a subset $E'\subset [0,\infty)$ with $\int_{E'}dr<\infty$ in the case $R_0=\infty$.
Thus, from (\ref{3.13}) and (\ref{3.16}) we get the desired inequality of Theorem \ref{3.14}.
$$ (q-\Delta_{\mathcal D,V}(k+1)-\epsilon) T_f(r,r_0)\le\sum_{i=1}^{q}\dfrac{1}{d}N^{[L-1]}_{Q_i(\tilde f)}(r)+S(r). $$

Return to the case (b). From Theorem \ref{3.14}, we have
$$\sum_{j=1}^q\delta^{[L-1]}_{f}(D_j)\le\sum_{j=1}^q\delta^{[L-1]}_{f,*}(D_j)\le \Delta_{\mathcal D,V}(k+1)+\epsilon.$$
The theorem is proved in this case.
\hfill$\square$

\noindent
{\bf Disclosure statement.} No potential conflict of interest was reported by the author(s).

\noindent
{\bf Acknowledgements.} We would like to thank the anonymous referee for a careful reading of the paper and for many helpful suggestions that helped us improve the presentation of the paper.

\vskip0.2cm
{\footnotesize 
\noindent
{\sc Tran Duc Ngoc}\\
Faculty of Mathematics and Applications, Saigon University,\\ 
 273 An Duong Vuong, Ward 3, Dist.5, Ho Chi Minh City, Vietnam.\\
\textit{E-mail}: tdngoc@moet.edu.vn}

\vskip0.2cm
{\footnotesize 
\noindent
{\sc Si Duc Quang}\\
Department of Mathematics, Hanoi National University of Education,\\
 136-Xuan Thuy, Cau Giay, Hanoi, Vietnam.\\
\textit{E-mail}: quangsd@hnue.edu.vn}

\end{document}